\documentclass[12pt,leqno,a4paper]{article}
\usepackage[sumlimits,namelimits]{amsmath}
\usepackage{amsthm,amssymb,amstext, amsmath}
\usepackage[T1]{fontenc}
\usepackage[latin1]{inputenc}
\usepackage{graphicx}
\usepackage{a4wide}
\usepackage[english]{babel}
\usepackage{amsfonts}
\usepackage{amsmath}
\usepackage{amssymb}
\usepackage{dsfont}

\newtheoremstyle{myplain}{6pt}{6pt}{\itshape}{}{\bfseries}{:}{.5em}{}
\newtheoremstyle{mydef}{6pt}{6pt}{}{}{\bfseries}{:}{\newline}{}
\newtheoremstyle{myrem}{6pt}{6pt}{}{}{\bfseries}{:}{.5em}{}

\theoremstyle{myplain}

\theoremstyle{mydef}

\theoremstyle{myrem}

\allowdisplaybreaks[3]
\newcommand {\e} {\varepsilon}

\newcommand {\ve} {\varepsilon}

\normalsize

\title{\bf On convergence and growth of sums $\sum c_k f(kx)$}
\author{Istvan Berkes\footnote{Graz University of Technology,
Institute of Statistics, Kopernikusgasse 24, 8010 Graz, Austria. \ \
\mbox{e-mail}: \texttt{berkes@tugraz.at}. Research supported by Austrian Science Fund (FWF), Grant
P24302. \newline {\it 2010 Mathematics Subject Classification.} 42A20, 42A61, 11K38.}
}
\begin{document}

\date{}
\maketitle

\abstract{We investigate the almost everywhere  convergence of  $\sum_{k=1}^\infty  c_k f(kx)$ where $f$ is a mean zero periodic function with bounded variation. The classical approach,
going back to the 1940's, depends on estimates for GCD sums $\sum_{k, \ell=1}^N (n_k, n_\ell)^2/(n_k n_\ell)$,  a connection leading to the a.e.\ convergence condition $\sum_{k=1}^\infty c_k^2 (\log k)^{2+\e}<\infty$.  Using a recent profound GCD estimate, Aistleitner and Seip \cite{aise} weakened the convergence condition to $\sum_{k=1}^\infty c_k^2 (\log k)^\e<\infty$, $\e>0$. In this paper we show that $\sum_{k=1}^\infty c_k^2 (\log \log k)^\gamma<\infty$ suffices for $\gamma> 4$ and not for $\gamma< 2$, settling the convergence problem except the unknown critical exponent of the loglog. Our method yields also new information on the growth of sums $\sum_{k=1}^N  f(n_kx)$, a problem closely connected with metric discrepancy estimates for $\{n_kx\}$. In analogy with the previous result we show that  optimal bounds for $\sum_{k=1}^N  f(n_kx)$ for $f \in BV$ differ from the analogous
(known) trigonometric results by a loglog factor.}

\section{Introduction}

Let $f:{\mathbb R}\to {\mathbb R}$ be a measurable function  satisfying
\begin{equation}\label{f0}
f(x+1)=f(x), \quad \int_0^1 f(x) dx=0, \quad \int_0^1 f^2(x) dx<\infty.
\end{equation}
The almost everywhere convergence of series
\begin{equation}\label{1}
\sum_{k=1}^\infty c_k f(kx)
\end{equation}
has been a much investigated problem of analysis. By Carleson's theorem \cite{ca}, in the case $f(x)=\sin 2\pi x$,
$f(x)=\cos 2\pi x$ the series (\ref{1}) converges a.e.\ provided $\sum_{k=1}^\infty c_k^2<\infty$. Gaposhkin \cite{ga68}
showed that this remains valid if the Fourier series of $f$ converges absolutely; in particular this holds
if $f$ belongs to the Lip ($\alpha$) class for some  $\alpha>1/2$. However, Nikishin \cite{ni} showed that the analogue
of Carleson's theorem fails for $f(x)=\hbox{sgn} \sin 2\pi x$ and fails also for some continuous $f$. Berkes \cite{be}
showed that there is a counterexample in the Lip (1/2) class, and thus Gaposhkin's result is sharp. For the class
Lip ($\alpha$), $0<\alpha\le 1/2$
and other classical function classes like $C(0, 1)$, $L^p(0, 1)$, $BV(0, 1)$ or function classes defined by the
order of magnitude of their Fourier coefficients, the results are much less complete: there are several
sufficient criteria (see e.g. \cite{ai}, \cite{bewe}, \cite{bewe2}, \cite{br}, \cite{ga66a}, \cite{we}; see also
\cite{ga66b} for the history of the subject until 1966) and necessary criteria (see \cite{be}), but there are large
gaps between the
sufficient and necessary conditions. Very recently, Aistleitner and Seip \cite{aise} proved that for $f\in \text{BV}$
the series (\ref{1}) converges a.e.\ if $\sum_{k=1}^\infty c_k^2 (\log k)^\e <\infty$ for some $\e>0$.
The purpose of this paper is to improve this result further and to provide a near optimal a.e.\ convergence condition
for $\sum_{k=1}^\infty  c_k f(kx)$ in the case $f\in BV$.

\medskip\noindent
{\bf Theorem 1.}  {\it Let $f: {\mathbb R}\to {\mathbb R}$ satisfy (\ref{f0}) and have bounded variation
on $[0, 1]$. Let $(c_k)$ be a real sequence satisfying
\begin{equation}\label{2}
\sum_{k=1}^\infty  c_k^2 (\log\log k)^\gamma<\infty
\end{equation}
for some $\gamma> 4$. Then for any increasing sequence $(n_k)$ of positive integers the series
$\sum_{k=1}^\infty c_k f(n_kx)$ converges a.e. On the other hand, for any $0<\gamma< 2$ there
exists an increasing sequence $(n_k)$ of positive integers and a real sequence $(c_k)$  such
that (\ref{2}) holds,
but $\sum_{k=1}^\infty c_k f(n_kx)$ is a.e.\ divergent for $f(x)=x-[x]-1/2$.}

\bigskip
The proof of Theorem 1 also provides crucial new information on the growth of sums
\begin{equation}\label{fsum}
\sum_{k=1}^N f(n_kx)
\end{equation}
where $f\in BV$ satisfies (\ref{f0}).  The order of magnitude of such sums is closely
related to the classical problem of estimating the discrepancy of $\{n_kx\}$ in the
theory of Diophantine approximation. The problem goes back to Weyl \cite{wey} and the
the strongest result for general $(n_k)$ is due to Baker \cite{ba}, who proved
that the discrepancy $D_N(\{n_kx\})$ of the sequence $\{n_kx\}_{1\le k \le N}$ satisfies
\begin{equation} \label {baker}
D_N (\{n_kx\}) \ll N^{-1/2} (\log N)^{3/2+\e} \quad \text{a.e.}
\end{equation}
for any $\e>0$ and any increasing sequence $(n_k)$ of integers. On the other hand,
Berkes and Philipp \cite{beph} constructed an increasing $(n_k)$ such that
$$ D_N (\{n_kx\}) \ll N^{-1/2} (\log N)^{1/2} \quad \text{a.e.} $$
is not valid. The optimal exponent of the logarithm remains open, in fact we do not even know
if (\ref{baker}) holds with $\e=0$. There has been, however, considerable progress in a closely
related, slightly easier problem. Relation (\ref{baker}) and Koksma's inequality
(see e.g. \cite{kn}, p.\ 143)
imply that for any $f\in BV $ satisfying (\ref{f0}) we have
\begin{equation*}
\left|\sum_{k=1}^N f(n_kx)\right| =O(\sqrt{N} (\log N)^{3/2+\e}) \quad \text{a.e.}
\end{equation*}
Aistleitner, Mayer and Ziegler \cite{aihozi} improved here the upper bound to
$$O(\sqrt{N} (\log N)^{3/2} (\log\log N)^{-1/2+\e})$$
and in an unpublished manuscript Berkes \cite{be3} showed that for polynomially increasing
$(n_k)$ the upper bound can be improved to $O(\sqrt{N} (\log N)^{1/2+\e})$. In \cite{aise} Aistleitner
and Seip removed the restriction of polynomial growth of $(n_k)$, obtaining the result for all $(n_k)$.
On the other hand, the sequence $(n_k)$ in Berkes and Philipp \cite{beph} actually satisfies
$$
\limsup_{N\to \infty} \frac{\left| \sum_{k=1}^N f(n_kx)\right|}{\sqrt{N} (\log N)^{1/2}}=\infty \qquad
\text{a.e.}
$$
for $f(x)=x-[x]-1/2$.
In this paper we will prove the following result.

\bigskip\noindent
{\bf Theorem 2.} {\it Let $f$ satisfy (\ref{f0}) and have bounded variation on $[0, 1]$. Let $\varphi$
be a  non-decreasing function satisfying $\sup_{k\ge 1} \varphi(2k)/\varphi (k)<\infty$ and
\begin{equation}\label{bp}
\sum_{k=1}^\infty \frac{1}{k\varphi(k)^2}<\infty.
\end{equation}
 Then for any increasing sequence $(n_k)$ of positive integers we have
\begin{equation}\label{th2}
\left|\sum_{k=1}^N f(n_kx)\right| =O(\sqrt{N} \varphi (N) (\log\log N)^2) \qquad \text{a.e.}
\end{equation}
}

\bigskip
To clarify the meaning of this theorem, let us note that Carleson's theorem combined with the
Kronecker lemma yields that under (\ref{bp}) we have
$$\lim_{N\to \infty} \frac{\sum_{k=1}^N \cos 2\pi n_k x}
{\sqrt{N}\varphi (N)}=0 \qquad \text{a.e.}
$$
Berkes and Philipp \cite{beph}  proved that this result is best possible in the sense that
if $\varphi$ is a non-decreasing function with
$\sup_{k\ge 1}\varphi (k^2)/\varphi (k)<\infty$
and
$$\sum^\infty_{k=1}\frac{1}{k\varphi (k)^2}=\infty,$$
then there exists an increasing sequence $(n_k)$ of integers such that
$$
\limsup_{N\to\infty}\frac{\left|\sum_{k=1}^N \cos2\pi n_kx \right|}
{\sqrt{N} \varphi (N)}=\infty\qquad\text{a.e.}
$$
This result describes precisely the growth of sums (\ref{fsum}) in the trigonometric case
$f(x)=\cos 2\pi x$.  Theorem 2 shows that, in analogy with Theorem 1, optimal bounds for
sums (\ref{fsum}) in the trigonometric case and for $f\in BV$ differ only in a loglog power.

\section{Proofs}

\medskip\smallskip
For the proof of Theorem 1 we need the following variant of Lemma 2 in \cite{be}.

\bigskip
\noindent{\bf Lemma.} {\it
Let $1\le p_1<q_1<p_2<q_2<\dots$ be integers such that $p_{k+1}\ge 4q_k$; let
$I_1, I_2,\dots$ be sets of integers such that $I_k\subset [2^{p_k}, 2^{q_k}]$
and each element of $I_k$ is divisible by $2^{p_k}$. Let $f(x)=x- [x]-1/2$ and
$$
X_k= X_k(\omega)=\sum_{j\in I_k} f(j\omega)\qquad
(k=1,2,\dots ,\,\,\,  \omega \in (0, 1) ).
$$
Then there exist independent r.v.'s $Y_1,Y_2,\dots$ on the probability space
$((0,1),\cal B, \lambda)$ such that $|Y_k|\le \text{card}\, I_k$, $EY_k=0$ and
$$
\|X_k-Y_k\|\le 2^{-k}\qquad (k\ge k_0),
$$
where $\| \cdot \|$ denotes the $L^2(0, 1)$ norm.
}\medskip

\medskip
\noindent {\bf Proof.}
Let ${\cal F}_k$ denote the $\sigma$-field generated by the dyadic intervals
\begin{equation}\label{4}
U_\nu=\left[ \nu 2^{-4q_k},  (\nu +1)2^{-4q_k}\right]\qquad
0\le \nu < 2^{4q_k}
\end{equation}
and set
\begin{align*}
& \xi_j = \xi_j(\cdot) =E(f(j\cdot)|{\cal F}_k), \qquad j\in I_k \\
& Y_k=Y_k(\omega)  =\sum_{j\in I_k} \xi_j(\omega).
\end{align*}
Since $|f|\le 1$, we have $|\xi_j|\le 1$ and thus $|Y_k|\le \text{card}\, I_k$.
Further, by $f\in \text{BV}$ the Fourier coefficients of $f$ are $O(1/k)$ and thus from Lemma 3.1 of \cite{be76} it follows that
$$
\|\xi_j (\cdot) -f(j\cdot)\| \le
C_1(j2^{-4q_k})^{1/3}\qquad j\in I_k
$$
and since $I_k$ has at most $2^{q_k}$ elements, we get
$$
\|X_k-Y_k\|\le C_2 \cdot 2^{-{q_k}} \le 2^{-k}
\qquad\text{for}\qquad
k \ge k_0.
$$
Since $p_{k+1}\ge 4q_k$ and since each $j\in I_{k+1}$ is a multiple
of $2^{p_{k+1}}$, each interval $U_\nu$ in (\ref{4}) is a period interval for
all $f(jx)$, $j\in I_{k+1}$ and thus also for $\xi_j$, $j\in I_{k+1}$.
Hence $Y_{k+1}$ is independent of the $\sigma$-field ${\cal F}_k$ and
since ${\cal F}_1\subset {\cal F}_2 \subset \dots$ and $Y_k$ is ${\cal F}_k$
measurable, the r.v.'s $Y_1, Y_2, \dots$ are independent.
Finally $E\xi_j=0$ and thus $EY_k=0$.
\medskip

\medskip
\noindent{\bf Proof of Theorem 1.}
We start with the second statement and actually prove a little more: we show that for any positive sequence $\e_k \to 0$
there exists an increasing sequence $(n_k)$ of integers and a real sequence $(c_k)$ such that
$$ \sum_{k=1}^\infty c_k^2 (\log\log k)^2 \ve_k<\infty$$
and $\sum_{k=1} ^\infty c_k f(n_k x)$ diverges a.e. for $f(x)=x-[x]-1/2$. Put $\e_k^*=\sup_{j\ge k} \e_j$ and
let $\psi(k)$ be a sequence of positive integers
growing so rapidly that $\psi(k+1)/\psi(k)\ge 2$ $(k=1, 2, \ldots)$ and setting
$$r_k=\psi(k)^3, \qquad M_k = \sum_{j \leq k} r_j \psi (j) = \sum_{j \leq k} \psi (j)^4, $$
we have
$$\sum_{k=1}^\infty \e_{M_{k-1}}^*<\infty.$$
By a well known result of G\'al \cite{gal}
there exists, for each $k\ge 1$, a sequence $m_1^{(k)} < m_2^{(k)} <\ldots < m_{\psi(k)}^{(k)}$
of positive integers such that
\begin{equation}\label{7}
\int\limits_0^1
\left( \sum_{j=1}^{\psi(k)}  f(m_j^{(k)} x) \right)^2 dx \ge
\text{const} \cdot \psi(k) (\log\log \psi(k))^2.
\end{equation}
We define sets
\begin{equation}\label{8}
I_1^{(1)},I_2^{(1)},\dots,I_{r_1}^{(1)},
I_1^{(2)},\dots,I_{r_2}^{(2)},
\dots,
I_1^{(k)},\dots,I_{r_k}^{(k)},
\dots
\end{equation}
of positive integers by
$$
I_j^{(k)}= 2^{c_j^{(k)}}
\left \{m_1^{(k)},\dots, m_{\psi(k)}^{(k)} \right\},\qquad
1\le j\le r_k,\; k\ge 1
$$
where $c_j^{(k)}$ are suitable positive integers. (Here for any set $\{ a,b,\dots\}\subset
R$ and $\lambda\in R$, $\lambda \{ a,b,\dots\}$ denotes the set
$\{ \lambda a,\lambda b,\dots\}$.) Clearly we can choose the integers
$c_j^{(k)}$ inductively so that the sets $I_j^{(k)}$ in (\ref{8}) satisfy the conditions assumed in
the Lemma for the sets $I_k$. Since the left hand side of (\ref{7}) does not change if we replace
every $m_j^{(k)}$ with $am_j^{(k)}$ for some integer $a \geq 1$, setting
$$
X_j^{(k)} = X_j^{(k)} (\omega)=
\sum_{\nu\in I_j^{(k)}} f (\nu\omega)
$$
we have
\begin{equation}\label{9}
E\left( X_j^{(k)}\right)^2 \ge \text{const}\cdot \psi(k) (\log\log \psi(k))^2.
\end{equation}
By the Lemma, there exist independent r.v.'s $Y_j^{(k)}$ ($1\le j\le r_k$,
$k=1,2,\dots$) such that $|Y_j^{(k)}|\le \psi(k)$, $EY_j^{(k)}=0$ and
\begin{equation}\label{10}
\sum_{k,j} \|X_j^{(k)}-Y_j^{(k)}\|<\infty \qquad
\end{equation}
whence
\begin{equation}\label{11}
\sum_{k,j} |X_j^{(k)}-Y_j^{(k)}|<\infty \quad \text{a.e.}
\end{equation}
By (\ref{9}) and (\ref{10}) we have
\begin{equation}\label{12}
E(Y_j^{(k)})^2 \ge \text{const} \cdot \psi(k) (\log\log \psi (k))^2.
\end{equation}
Hence setting
$$
Z_k = \dfrac{1}{\sqrt{r_k\psi(k)} \log\log \psi (k)} \sum_{j=1}^{r_k} Y_j^{(k)}, \qquad
\sigma_k^2 = E\left( \sum_{j=1}^{r_k} Y_j^{(k)}\right)^2,
$$
we get from the central limit theorem with Berry-Esseen remainder term (see e.g. Feller \cite{fe}, p.\ 544),
(7), and $r_k=\psi(k)^3$,
$$
\gathered
P(Z_k \ge 1)\ge P\left( \sum_{j=1}^{r_k} Y_j^{(k)} \ge c \sigma_k \right)
\ge\\
\ge (1-\Phi (c))-\text{const}\cdot \dfrac{r_k \psi(k)^3}{(r_k\psi(k)
(\log\log \psi(k))^2)^{3/2}}\ge\\
\ge 1-\Phi (c)- o(1)\ge c'>0\qquad (k\ge k_0)
\endgathered
$$
where $\Phi$ denotes the Gaussian distribution function and $c$ and $c'$
are positive absolute constants. Since the r.v.'s $Z_k$ are
independent, the Borel--Cantelli lemma implies that $P(Z_k\ge 1 \; \text{i.o.})=1,$
i.e. $\sum_{k\ge 1} Z_k$ is a.e.\ divergent, which, in view of (\ref{11}), yields that
\begin{equation}\label{13}
\sum_{k=1}^\infty\,
\dfrac{1}{\sqrt{r_k \psi(k)} \log\log \psi(k)}\,
\sum_{j=1}^{r_k} X_j^{(k)}
\qquad\text{is a.e.\ divergent.}
\end{equation}
In other words, $\sum_{i=1}^\infty c_i f(n_i x)$ is a.e.\ divergent, where
\begin{align}\label{14}
&(n_i)_{i\ge 1}= \bigcup\limits_{k=1}^\infty\,
\bigcup\limits_{j=1}^{r_k} I_j^{(k)}
\end{align}
and
$$ c_i^2 =  \frac{1}{r_k \psi (k)
(\log\log \psi (k))^2} \qquad \text{for} \quad M_{k-1} < i \leq M_k.$$
Now for $M_{k-1} < i \leq M_k$ we have by the exponential growth of $\psi(j)$,
$$i \le 2\psi (k) ^4, \quad \log\log i \le
2\log\log \psi (k) \qquad (k\ge k_0)$$
and consequently for $M_{k-1} < i \leq M_k$
$$ c_i^2 (\log\log i)^2 \e_i \leq \text{const}
\cdot
\frac{1}{r_k \psi (k)} \e_{M_{k-1}}^*. $$
Hence
$$ \sum_{i=1}^{\infty} c_i^2 (\log\log i)^2 \varepsilon_i \leq \text{const} \cdot \sum_{k=1}^{\infty}
\e_{M_{k-1}}^*  < +\infty,$$
completing the proof of the second half of Theorem 1.

We prove now the sufficiency of (\ref{2}) in Theorem 1 for $\gamma>4$.
In [2], a slightly weaker result is proved, namely the convergence of
$\sum_{k=1}^\infty  c_k f(kx)$ under the assumption $\sum_{k=1}^\infty c_k^2 (\log k)^\e<\infty$
for some $\e>0$. To get the present result, a slight improvement of the argument in
\cite{aise} is needed. Let
$$\Gamma (N)=\sup_{n_1, \ldots, n_N} \frac{1}{N} \sum_{k, \ell=1}^N
\frac{(\text{gcd} (n_k, n_\ell))^{2\alpha}}{(n_k n_\ell )^\alpha},$$
where the supremum is taken over all distinct positive integers $n_1, \ldots, n_N$.
We start out of the formula
\begin{equation}\label{seip}
\Gamma (N)\le \prod_{j=1}^{r_N}(1-v_j)^{-1}(1-v_j^{-1}\tau_j^2)^{-1}
\prod_{k=r_N+1}^{N-1} (1-v_{r_N}^{-1}\tau_k^2)^{-1}+\exp \left( C\sum_{\ell=1}^{N-1} t_\ell^2\right)
\end{equation}
on page 7 of \cite{aise}. We have chosen here $\xi=2$, $C$ is an absolute constant  and letting $p_j$ denote the $j$-th prime, for $\alpha\in (\log 2/\log 3, 1)=I$ we have $p_j^{-\alpha}< 1/2$ for $j\ge 2$ and thus with the notations of \cite{aise} we have for $j\ge 2$
$$t_j=p_j^{-\alpha}, \quad \tau_j=2p_j^{-\alpha},  \quad v_j=\max ( \tau_j, (2\alpha-1)^{-1/2}
\tau_{r_N}),  \quad r_N=[2\log N] +1.$$
We estimate the right hand side of (\ref{seip}), just as in \cite{aise}, by using the prime number theorem.
Let
$$s_N=\max \{1\le j\le r_N: \tau_j\ge (2\alpha-1)^{-1/2} \tau_{r_N}\}$$
for $N\ge N_0$ and split the first product on the right hand side of (\ref{seip}) into two sub-products $P_1$ and $P_2$,
extended for the indices $1\le j \le s_N$ and $s_N+1 \le j \le r_N$. Let further $P_3$ denote the second product
on the right hand side of (\ref{seip}). We estimate $P_1, P_2, P_3$ separately. Letting $C_1, C_2, \ldots$
denote positive absolute constants, let us observe that by $p_k\sim k\log k$ we have for any $\alpha\in I$
\begin{align*}
&\prod_{j=4}^{s_N}(1-2p_j^{-\alpha})^{-2}\le \prod_{j=4}^{r_N}(1-2p_j^{-\alpha})^{-2}
\le \exp\left( 4\sum_{j=4}^{r_N} p_j^{-\alpha}\right)\le \exp\left( C_1\sum_{j=4}^{r_N} (j\log j)^{-\alpha}\right)\\
& \le \exp\left( C_1\int_3^{r_N} (x\log x)^{-\alpha}\, dx\right)\le \exp\left( C_1\int_3^{r_N} x^{-\alpha}\, dx\right)\le
\exp\left( \frac{C_1}{1-\alpha} r_N^{1-\alpha}\right)\\
&\le \exp\left( \frac{C_2}{1-\alpha} (\log N)^{1-\alpha}\right).
\end{align*}
Thus we have
\begin{equation*}
P_1 \le C_3 \exp \left( \frac{C_2}{1-\alpha} (\log N)^{1-\alpha}\right)
\quad \text{for} \ \alpha\in I.
\end{equation*}
Similar estimates hold for $P_2$ and $P_3$ (which do not involve a singularity at $\alpha=1$ and hence the factor
$1/(1-\alpha)$) and we arrive at
\begin{equation}\label{gamma}
\Gamma (N) \le \exp \left( \frac{C_4}{1-\alpha} (\log N)^{1-\alpha}\right) \quad  \text{for}
\ \alpha \in I.
\end{equation}
We thus see that if in the first estimate of Theorem 1 of \cite{aise} we drop the factor $(\log\log N)^{-\alpha}$,
the resulting estimate holds uniformly for $\alpha \in I$.

In relation (26) of \cite{aise} the arbitrary parameter
$0<\e<1$ appears. The subsequent
arguments lead to relation (27) in \cite{aise}, yielding the norm bound
\begin{equation}\label{5}
cJ^{-\e/2} (\log N)^{1/2}\left( \sum_{k=M_1+1}^{M_2} c_k^2\right)^{1/2} \exp (c(\e)
(\log N)^{\e/2}(\log\log N)^{-1/2}).
\end{equation}
Using (\ref{gamma}) instead of the estimate in the first line of Theorem 1 in \cite{aise},
we get
\begin{equation}\label{6}
C_5J^{-\e/2} (\log N)^{1/2}\left( \sum_{k=M_1+1}^{M_2} c_k^2\right)^{1/2} \exp \left(\frac{C_6}{\e}
(\log N)^{\e/2}\right)
\end{equation}
instead of (\ref{5}). By increasing $C_6$ if necessary, we can assume $C_6\ge 4$. Here $C_5$ and the further constants
$C_7, C_8, \ldots$ are allowed to depend also on $f$. We choose now $J$ by
$$J^{\e/2}= (\log N)^{1/2}\exp \left(\frac {2C_6}{\e}(\log N)^{\e/2}\right)$$
and thus the expression (\ref{6}) becomes
$$C_5\left( \sum_{k=M_1+1}^{M_2} c_k^2\right)^{1/2} \exp \left(-\frac {C_6}{\e}(\log N)^{\e/2}\right).$$
Applying the Rademacher-Mensov inequality as in \cite{aise}, it follows that the norm in formula (29) of
\cite{aise}
can be bounded by
\begin{equation}\label{rm}
C_7 (\log N) \exp \left(-\frac {C_6}{\e}(\log N)^{\e/2}\right) \left( \sum_{k=1}^N c_k^2\right)^{1/2}.
\end{equation}
Choosing $\e=1/(\log\log N)$ and using $C_6\ge 4$,
the expression in (\ref{rm}) will be $\le C_8 (\sum_{k=1}^N c_k^2)^{1/2}$. On the other hand,
\begin{equation}\label{ca}
\log J=\frac{1}{\e} \log\log N+ \frac{4C_6}{\e^2} (\log N)^{\e/2},
\end{equation}
which becomes $\le C_9 (\log\log N)^2$ with the choice of $\e$. Thus the last expression in the first
displayed formula on page 16 of \cite{aise} becomes
$$\le C_{10} (\log\log N)^2 \left( \sum_{k=1}^{N} c_k^2\right)^{1/2}.$$
This shows that the expression $c(\e) (\log N)^\e$ on the right hand side
of the maximal inequality in Lemma 4 of \cite{aise} can be replaced by $C_{11} (\log\log N)^4$.
Arguing further as in \cite{aise}, this leads to the sufficiency of the convergence condition (\ref{2}) for $\gamma>4$.

\bigskip\noindent
{\bf Proof of Theorem 2.} Let $S_N=S_N(x)= \sum_{k=1}^N f(n_kx)$ and let $\varphi$ be a non-decreasing function
satisfying (\ref{bp}) and $\sup_{k\ge 1} \varphi(2k)/\varphi (k)<\infty$.  By the just proved stronger
version of Lemma 4 of \cite{aise} we have
\begin{equation*}
\int_0^1 \max_{1\le k \le N} S_k^2\, dx \le C N (\log\log N)^4  \end{equation*}
with some constant $C>0$. Thus
$$P\left(\max_{1\le k \le 2^N} |S_k|\ge \sqrt{2^N} \varphi (2^N) (\log\log 2^N)^2\right)
\le \frac{C}{\varphi^2 (2^N)}.$$
Now (\ref{bp}) and the monotonicity of $\varphi$ imply
$\sum_{k=1}^\infty \frac{1}{\varphi^2 (2^k)}<\infty$ and thus Theorem 2 follows
from the Borel-Cantelli lemma.

\end{document}